# APPLICATION OF TRIPLE SHEHU TRANSFORMS TO FRACTIONAL DIFFERENTIAL EQUATIONS

WAGDI F. S. AHMED AND D. D. PAWAR

ABSTRACT. The triple Shehu transform, a new generalisation of the triple Laplace transform and triple Sumudu transform, has recently been introduced. The triple Shehu transform formulas for fractional Caputo operators were obtained in this study. The generalised integral transform was subsequently applied to solve fractional partial differential equations including the Caputo derivative.

## 1. INTRODUCTION

In the fields of physics, engineering, and other applied disciplines, the application of differential equations is one of the most attractive and important topics. There are no standard approaches for solving such issues involving differential equations. The integral transform technique is oneof the most well-known schemes that is used by many scholars for solving ordinary and partial differential equations [3, 4, 5, 6, 7, 16, 24, 29, 31]. Several approaches have been introduced and developed to get approximate or exact solutions of these types of equations including, the Adomian de-composition method (ADM)[1, 2], fractional residual power series method [33], variational iteration method[44], the homotopy perturbation method[39] differential transform method[30]. Among these methods, the integral transformation methods which are rather and popular, thus in the literature, there are various types of integral transforms such as Laplace transform [11, 13, 21, 22, 32, 35, 36, 41], Mellin transform[40], Sumudu transform[25, 42, 43], and many others. These transforms are widely implemented toget analytical solutions of ordinary, partial differential equations, and integral equations. After the appearance of one dimensional integral transform,the idea was extended to the two dimensional and named the double integral transform [23, 26, 28]. One dimensional Shehu transform [38] is closely related to Laplace and Sumudu transforms. Shehu Maitama and Weidong Zhao introduced this transform for the first time in 2019 and used it to solve ordinary and partial differential equations. For additional information on the Shehu transform, its characteristics, and its uses, see [18, 20, 34]. Double Shehu transform [12] which is a new generalization of double Sumudu trans-form and double Laplace transform, was recently introduced by Alfaqeih and MĂśsĂśrlĂś which applied to get the exact solutions of partial differential equations of two variables and integral equations with convolution.

Triple Shehu transform [10] which is a new generalization of triple Laplace transform, was recently introduced by Alkaleeli et al which applied to get the exact solutions of fractional partial differential equations of three variables. The organization of the article is presented as follows: The second section, we give some notations about fractional calculus. In section 3, we present the definition of one dimensional Shehu transform, the double Shehu trans- form and the triple Shehu transform, with related theorems. In section 4, we derive the formulas of the triple Shehu transform for the fractional Caputo operators. In section 5, we implement the triple Shehu transform to solve fractional differential equations. The conclusion is presented in section 6.

## 2. PRELIMInaRIES

This section presents some fundamental notations for Shehu transformsand fractional calculus.

**Definition 2.1.**[14,15] Let $f(x,t)$ be a continuous function. Then,the Riemann-Liouville (RL) partial fractional integrals are given by :

$$ {}_tI^{\gamma} f(x,t) = \frac{1}{\Gamma(\gamma)} \int_0^t (t-\tau)^{\gamma-1} f(x,\tau)d\tau, \quad (1)$$

$$ {}_xI^{\beta} f(x,t) = \frac{1}{\Gamma(\beta)} \int_0^x (x-\lambda)^{\beta-1} f(\lambda,\tau)d\lambda, \quad (2)$$

Where $(x,t) \in (0,\infty) * (0,\infty)$ $and$ $\beta, \gamma > -1$

**Definition 2.2** [15] Let f(x, t) be a continuous function. Then, the Riemann-Liouville(RL) partial fractional derivatives are given by:

$$ {}_RD_t^{\gamma} f(x,t) = \frac{1}{\Gamma(r-\gamma)} (\frac{d}{dt})^r \int_0^t (t-\tau)^{r-\gamma-1} f(x,\tau)d\tau, \quad (3)$$

$$ {}_RD_x^{\beta} f(x,t) = \frac{1}{\Gamma(r-\beta)} (\frac{d}{dt})^n \int_0^x (x-\lambda)^{n-\beta-1} f(\lambda,t)d\lambda, \quad (4)$$

**Definition 2.3** [15] The Caputo partial fractional derivatives of a function f(x, t) are defined by:

$$ D_t^{\gamma} f(x,t) = \frac{1}{\Gamma(r-\gamma)} \int_0^t (t-\tau)^{r-\gamma-1} f(x,\tau)d\tau, \quad (5)$$

$$ D_x^{\beta} f(x,t) = \frac{1}{\Gamma(r-\beta)} \int_0^x (x-\lambda)^{n-\beta-1} f(\lambda,t)d\lambda, \quad (6)$$

Where $r-1 < \gamma < r, n-1 < \beta < n, and\ r,n \in N$

The following are some basic properties of the Caputo derivatives and partialfractional integrals

(i) $ {}_tI^{\gamma}\ {}_xI^{\beta} f(x,t) = \frac{1}{\Gamma(\gamma)\Gamma(\beta)} \int_0^t (t-\tau)^{\gamma-1}(x-\lambda)^{\beta-1} f(\lambda,\tau)d\tau d\lambda,$

(ii) $D_t^{\gamma} D_x^{\beta} f(x,t) = \frac{1}{\Gamma(r-\beta)\Gamma(r-\gamma)} \int_0^t (t-\tau)^{r-\gamma-1}(x-\lambda)^{n-\beta-1} f^{(r+n)}(\lambda,\tau)d\tau d\lambda,$

(iii) $ {}_tI^{\gamma} D_t^{\gamma} f(x,t) = f(x,t) - \sum_{i=0}^{n} \frac{f^{(i)}(x,0)}{i!} t^i,$

(iv) $ {}_tI^{\gamma} D_t^{\gamma} f(x,t) = f(x,t)$

**Definition 2.4.** The Mittag – Leffler function $E_{\gamma,\beta}(z)$ was defined as $E_{\gamma,\beta}(x) = \sum_{r=0}^{\infty} \frac{x^r}{\Gamma(\gamma r+\beta)} \quad (\beta, \gamma, x \in C, R(\gamma) > 0, R(\beta) > 0)$ (7)

**Definition 2.5.** In 1961, Fox defined the Ή -function as the Mellin –Barnes type path integral

$$Ή_{u,v}^{n,m}\left[-\sigma \left| \begin{matrix} (a_\gamma, A_\gamma)_1^u \\ (b_\gamma, B_\gamma)_1^v \end{matrix} \right.\right] = \frac{1}{2\pi i}\int_L \frac{\prod_{\gamma=1}^{n}(b_\gamma - B_\gamma\tau)\prod_{\beta=1}^{m}\Gamma(1-\alpha_\beta+\tau A_\beta)\sigma^\tau d\tau}{\prod_{\gamma=n+1}^{v}\Gamma(1-b_\gamma+B_\gamma\tau)\prod_{\beta=m+1}^{u}(\alpha_\beta-\tau A_\beta)}, \quad (8)$$

Where L is a suitable contour the orders $(n, m, u, v)$ are integers $0nv$, $0mu$, and the parameters $a_\beta \in R, A_\beta > 0, \beta = 1,2,3, \dots, u, b_\gamma > 0, b_\gamma \in R_\gamma = 1,2,3, \dots, v$ are such that

$$A_\beta(b_\gamma + \alpha) \neq B_\gamma(\alpha_\beta - \alpha - 1), \alpha - 1,2,3, \dots$$

$$Ή_{\tau,z+1}^{1,\tau}\left[-\sigma \left| \begin{matrix} (1-a1, A_1), \dots, (1-a\tau, A_\tau) \\ (1,0), (1-b1, B_1), \dots, (1-bz, B_z) \end{matrix} \right.\right]$$

$$= \sum_{s=0}^{\infty} \frac{\Gamma(a_1+A_1 s), \dots, \Gamma(a_\tau+A_\tau s)}{s!\Gamma(b_1+B_1 s), \dots, \Gamma(b_z+B_z s)} \quad (9)$$

## 3. SHEHU TRANSFORM

**Definition 3.1.** [37]The single Shehu transforms (HT) of a real-valued func- tion f (x, y, t) with respect to the variables x,y and t respectively, are defined by:

$$H_x(f(x,y,t)) = \int_0^\infty e^{-\left[\frac{ax}{h}\right]} f(x,y,t)dx, \quad (10)$$

$$H_y(f(x,y,t)) = \int_0^\infty e^{-\left[\frac{by}{k}\right]} f(x,y,t)dy, \quad (11)$$

$$H_t(f(x,y,t)) = \int_0^\infty e^{-\left[\frac{ct}{l}\right]} f(x,y,t)dt, \quad (12)$$

**Definition 3.2.**[8] The double Shehu transform of the function $f(x,t)$ isgiven by

$$H_x H_y(f(x,y)) = F[(a,b),(h,k)] = \int_0^\infty \int_0^\infty e^{-\left[\frac{ax}{h}+\frac{by}{k}\right]} f(x,y)dxdy, \quad (13)$$

**Lemma 3.3.**[19] Let $\boldsymbol{\gamma > 0}$, ***and*** $\boldsymbol{f(x,y)}$ are of exponential order. Then, thesingle Shehu transform of ${}_y\mathbf{I}^\gamma \mathbf{f(x,y)}$ was defined as:

$$H_y({}_yI^\gamma f(x,y)) = \left(\frac{b}{k}\right)^{-\gamma} H_y f(x,y) \quad (14)$$

**Lemma 3.4.**[19] Let $\boldsymbol{\gamma > 0, n-1 < \gamma < n, (n \in N)}$ be such that $\boldsymbol{f \in C^r(0,\infty)}$ and is of exponential order. Then, the single Shehu transform ofCaputo fractional derivative $D^\gamma f(x,y)$ was defined as:

$$H_y\left(D_y^{\gamma} f(x,y)\right)=\left(\frac{b}{k}\right)^{\gamma} H_y f(x,y)-\sum_{i=0}^{n-1}\left(\frac{b}{k}\right)^{\gamma-1-i}\left(\frac{\partial^i}{\partial y^i} f(x,0)\right) \tag{15}$$

**Theorem 3.5.** [19] The single Shehu transform of $y^{\beta-1}E_{\gamma-\beta}(cy^{\gamma})$ is givenby:

$$H_y\left(y^{\beta-1}E_{\gamma-\beta}(cy^{\gamma})\right)=\frac{\left(\frac{b}{k}\right)^{\gamma-\beta}}{\left(\frac{b}{k}\right)^{\gamma}-c},\ [c]<\left|\left(\frac{b}{k}\right)^{\gamma}\right| \tag{16}$$

**Theorem 3.6.** [9] Let $\boldsymbol{\gamma},\boldsymbol{\beta}>\mathbf{0}$ and $\boldsymbol{f}(\boldsymbol{x},\boldsymbol{y})$ is of exponential order. Then,the double Shehu transform of the fractional integrals $_{x}\boldsymbol{I}^{\boldsymbol{\gamma}}\ \boldsymbol{f}(\boldsymbol{x},\boldsymbol{y})$, $_{y}\boldsymbol{I}^{\boldsymbol{\beta}}\ \boldsymbol{f}(\boldsymbol{x},\boldsymbol{y})$ and $_{y}\boldsymbol{I}^{\boldsymbol{\beta}}\ \boldsymbol{f}(\boldsymbol{x},\boldsymbol{y})\ _{x}\boldsymbol{I}^{\boldsymbol{\gamma}}\ \boldsymbol{f}(\boldsymbol{x},\boldsymbol{y})$ is given by the follows:

$$H_xH_y\left(\ _{y}I^{\gamma}\ f(x,y)\right)=\left(\frac{b}{k}\right)^{-\gamma} H_y f(x,y) \tag{17}$$

$$H_xH_y\left(\ _{x}I^{\beta}\ f(x,y)\right)=\left(\frac{a}{h}\right)^{-\beta} H_y f(x,y) \tag{18}$$

$$H_xH_y\left(\ _{y}I^{\gamma}\ _{x}I^{\beta}\ f(x,y)\right)=\left(\frac{b}{k}\right)^{-\gamma}\left(\frac{a}{h}\right)^{-\beta} H_y f(x,y) \tag{19}$$

**Theorem 3.7.** [9] Let $\boldsymbol{\gamma},\boldsymbol{\beta}>\mathbf{0},\boldsymbol{n}-\mathbf{1}<\boldsymbol{\gamma}<\boldsymbol{n},\boldsymbol{r}-\mathbf{1}<\boldsymbol{\beta}<\boldsymbol{r}(\boldsymbol{r},\boldsymbol{n}\boldsymbol{\epsilon}\boldsymbol{N})$ be such that $\boldsymbol{f}\in\boldsymbol{C}^{\boldsymbol{l}}[(\mathbf{0},\infty)*(\mathbf{0},\infty)]$, order. Then, the double Shehu transforms of Caputo fractional derivativesare given by:

$$H_x\ H_y\left(D_y^{\gamma} f(x,y)\right)=\left(\frac{b}{k}\right)^{\gamma} H_xH_y f(x,y)-\sum_{i=0}^{n-1}\left(\frac{b}{k}\right)^{\gamma-1-i} H_x\left(\frac{\partial^i}{\partial y^i} f(x,0)\right) \tag{20}$$

$$H_x\ H_y\left(D_x^{\beta} f(x,y)\right)=\left(\frac{a}{h}\right)^{\beta} H_xH_y f(x,y)-\sum_{j=0}^{r-1}\left(\frac{a}{h}\right)^{\beta-1-j} H_y\left(\frac{\partial^j}{\partial y^j} f(0,y)\right) \tag{21}$$

$$H_x\ H_y\left(D_x^{\beta}D_y^{\gamma} f(x,y)\right)=\left(\frac{a}{h}\right)^{\beta}\left(\frac{b}{k}\right)^{\gamma} H_xH_y f(x,y)-\left(\frac{b}{k}\right)^{\gamma-1-i}\sum_{j=0}^{r-1}\left(\frac{a}{h}\right)^{\beta-1-j} H_y\left(\frac{\partial^j}{\partial y^j} f(0,y)\right)-\left(\frac{a}{h}\right)^{\beta}\sum_{i=0}^{n-1}\left(\frac{b}{k}\right)^{\gamma-1-i} H_x\left(\frac{\partial^i}{\partial y^i} f(x,0)\right)+\sum_{i=0}^{r-1}\ \sum_{j=0}^{r-1}\left(\frac{b}{k}\right)^{\gamma-1-i}\left(\frac{a}{h}\right)^{\beta-1-j}\left(\frac{\partial^{i+j}}{\partial x^j\partial y^i} f(0,0)\right) \tag{22}$$

**Definition 3.8.** [10] The Triple Shehu transform of the function $f(x,y,t)$ is given by

$$H_{x,y,t}^{-3}\left(f(x,y,t)\right)=F[(a,b,c),(h,k,l)]=\int_0^{\infty}\int_0^{\infty}\int_0^{\infty} e^{-\left[\frac{ax}{h}+\frac{by}{k}+\frac{ct}{l}\right]} f(x,y,t)\,dxdydt \tag{23}$$

Where $\mathbf{x},\mathbf{y},\mathbf{t}\geq\mathbf{0}$ and $\mathbf{a},\mathbf{b},\mathbf{c},\mathbf{h},\mathbf{k}$ **and** $\mathbf{l}$ are Shehu variables.

**Definition 3.9.** [10] The Triple inverse Shehu transform is given by

$$H_{x,y,t}^{-3}(F[(a,b,c),(h,k,l)])=f(x,y,t)$$

$$=\frac{1}{2\pi i}\int_{\gamma-i\infty}^{\gamma+i\infty}\frac{1}{h}e^{\left[\frac{ax}{h}\right]}\left[\frac{1}{2\pi i}\int_{\beta-i\infty}^{\beta+i\infty}\frac{1}{k}e^{\left[\frac{by}{k}\right]}\left(\frac{1}{2\pi i}\int_{\zeta-i\infty}^{\zeta+i\infty}\frac{1}{l}e^{\left[\frac{ct}{l}\right]}H_3(f(x,y,t))dc\right)db\right]da \tag{24}$$

**Theorem 3.10.**[10] The convolution of Triple Shehu transform of the function $f(x,y,t), g(x,y,t)$ is denoted by $(f *** g)(x,y,t)$ and defined by

$$(f *** g)(x,y,t) = \int_0^x \int_0^y \int_0^t f(x-z_1, y-z_2, t-z_3) g(z_1, z_2, z_3) dz_1 dz_2 dz_3$$

$$= \int_0^x \int_0^y \int_0^t g(x-z_1, y-z_2, t-z_3) f(z_1, z_2, z_3) dz_1 dz_2 dz_3 \qquad (25)$$

## 4. MAIN RESULTS

In this part, we establish the triple Shehu transform of the factional inte-grals and derivatives.

**Theorem 4.1.** Let $\gamma > 0$ and $f(x,y,t)$ is of exponential order. Then, the Triple Shehu transform of the fractional integrals ${}_tI^{\gamma} f(x,y,t)$ is given by the follows:

$$H_x H_y H_t\left( {}_tI^{\gamma} f(x,y,t)\right) = \left(\frac{c}{l}\right)^{-\gamma} H_x H_y H_t f(x,y,t) \qquad (26)$$

Proof. By using the Triple Shehu transform of the convolution with respectto t, we have

$$H_x H_y H_t\left( {}_tI^{\gamma} f(x,y,t)\right) = H_x H_y H_t\left(\frac{1}{\Gamma(\gamma)} t^{\gamma-1} ** f(x,y,t)\right)$$

$$= H_x H_y \left(\frac{1}{\Gamma(\gamma)} H_t(t^{\gamma-1}) H_t f(x,y,t)\right)$$

$$= H_x H_y \left(\frac{1}{\Gamma(\gamma)} \frac{(\gamma-1)!}{\left(\frac{c}{l}\right)^{\gamma}} H_t f(x,y,t)\right)$$

$$= \left(\frac{c}{l}\right)^{-\gamma} H_x H_y H_t f(x,y,t)$$

**Theorem 4.2.** Let $\boldsymbol{\beta} > \mathbf{0}$ and $\mathbf{f(x,y,t)}$ is of exponential order. Then, the Triple Shehu transform of the fractional integrals ${}_{\mathbf{y}}\mathbf{I}^{\boldsymbol{\beta}}\, \mathbf{f(x,y,t)}$ is given by the follows:

$$H_x H_y H_t\left( {}_yI^{\beta} f(x,y,t)\right) = \left(\frac{b}{k}\right)^{-\beta} H_x H_y H_t f(x,y,t) \qquad (27)$$

*Proof.* By using the Triple Shehu transform of the convolution with respectto y, we

have

$$H_x H_y H_t\left( {}_yI^{\beta} f(x,y,t)\right) = H_x H_y H_t\left(\frac{1}{\Gamma(\beta)} y^{\beta-1} ** f(x,y,t)\right)$$

$$= H_x H_t \left(\frac{1}{\Gamma(\beta)} H_y(y^{\beta-1}) H_y f(x,y,t)\right)$$

$$= H_x H_t \left(\frac{1}{\Gamma(\gamma)} \frac{(\beta-1)!}{\left(\frac{b}{k}\right)^{\beta}} H_y f(x,y,t)\right)$$

$$= \left(\frac{b}{k}\right)^{-\beta} H_x H_y H_t f(x,y,t)$$

**Theorem 4.3.** Let $\boldsymbol{\zeta > 0}$ and $\mathbf{f(x,y,t)}$ is of exponential order. Then, the Triple Shehu transform of the fractional integrals $\ {}_{\mathbf{x}}\mathbf{I}^{\boldsymbol{\zeta}}\, \mathbf{f(x,y,t)}$ is given by the follows:

$$H_x H_y H_t\left( {}_xI^{\zeta} f(x,y,t)\right) = \left(\frac{a}{h}\right)^{-\zeta} H_x H_y H_t f(x,y,t) \qquad (28)$$

*Proof.* By using the Triple Shehu transform of the convolution with respectto x, we have

$$H_x H_y H_t\left( {}_xI^{\zeta} f(x,y,t)\right) = H_x H_y H_t\left(\frac{1}{\Gamma(\zeta)} x^{\zeta-1} ** f(x,y,t)\right)$$

$$= H_y H_t\left(\frac{1}{\Gamma(\zeta)} H_x(x^{\zeta-1}) H_x f(x,y,t)\right)$$

$$= H_y H_t\left(\frac{1}{\Gamma(\zeta)} \frac{(\zeta-1)!}{\left(\frac{a}{h}\right)^{\zeta}} H_x f(x,y,t)\right)$$

$$= \left(\frac{a}{h}\right)^{-\zeta} H_x H_y H_t f(x,y,t)$$

**Theorem 4.4.** Let $\gamma, \beta, \zeta > 0$ and $f(x,y,t)$ is of exponential order. Then,theTriple Shehu transform of the fractional integrals $\ {}_tI^{\gamma}\ {}_yI^{\beta}\ {}_xI^{\zeta} f(x,y,t)$ are given by the follows:

$$H_x H_y H_t\left( {}_tI^{\gamma}\ {}_yI^{\beta}\ {}_xI^{\zeta} f(x,y,t)\right) = \left(\frac{c}{l}\right)^{-\gamma} \left(\frac{b}{k}\right)^{-\beta} \left(\frac{a}{h}\right)^{-\zeta} H_x H_y H_t f(x,y,t) \qquad (29)$$

Proof. By using the triple Shehu transform of the convolution we have

$$H_x H_y H_t\left( {}_tI^{\gamma}\ {}_yI^{\beta}\ {}_xI^{\zeta} f(x,y,t)\right) = H_x H_y H_t\left(\frac{1}{\Gamma(\gamma)} \frac{1}{\Gamma(\beta)} \frac{1}{\Gamma(\zeta)} x^{\zeta-1} y^{\beta-1} t^{\gamma-1} ** f(x,y,t)\right)$$

$$= \frac{1}{\Gamma(\gamma)} \frac{1}{\Gamma(\beta)} \frac{1}{\Gamma(\zeta)} H_x H_y H_t\left(H_x(x^{\zeta-1} y^{\beta-1} t^{\gamma-1}) H_x H_y H_t f(x,y,t)\right)$$

$$= \frac{1}{\Gamma(\gamma)} \frac{1}{\Gamma(\beta)} \frac{1}{\Gamma(\zeta)} H_x H_y H_t\left(\frac{(\zeta-1)!}{\left(\frac{a}{h}\right)^{\zeta}} \frac{(\beta-1)!}{\left(\frac{b}{k}\right)^{\beta}} \frac{(\gamma-1)!}{\left(\frac{c}{l}\right)^{\gamma}} H_y H_t H_x f(x,y,t)\right)$$

$$= \left(\frac{c}{l}\right)^{-\gamma} \left(\frac{b}{k}\right)^{-\beta} \left(\frac{a}{h}\right)^{-\zeta} H_x H_y H_t f(x,y,t)$$

**Theorem 4.5.** Let $\gamma, \beta, \zeta > 0, n-1 < \gamma < n, r-1 < \beta < r, m-1 < \zeta < m (r,n,m \epsilon N)$ be such that $f \in C^l[(0,\infty) * (0,\infty) * (0,\infty)]$, L= max{m,n,r} and is of exponential order. Then, the Triple Shehu transforms of Caputo fractional derivatives are given by:

$$H_t H_x\, H_y\left(D_x^{\zeta} f(x,y,t)\right) = \left(\frac{a}{h}\right)^{\zeta} H_t H_x\, H_y f(x,y,t) - \sum_{u=0}^{m-1} \left(\frac{a}{h}\right)^{\zeta-1-u} H_y H_t\left(\frac{\partial^u}{\partial y^u} f(0,y,t)\right) \quad (30)$$

$$H_t H_y\, H_x\left(D_y^{\beta} f(x,y,t)\right) = \left(\frac{b}{k}\right)^{\beta} H_t H_y H_x f(x,y,t) - \sum_{j=0}^{n-1} \left(\frac{b}{k}\right)^{\beta-1-j} H_t H_x\left(\frac{\partial^j}{\partial y^j} f(x,0,t)\right) \quad (31)$$

$$H_t H_y\, H_x\left(D_t^{\gamma} f(x,y,t)\right) = \left(\frac{c}{l}\right)^{\gamma} H_t H_y H_x f(x,y,t) - \sum_{i=0}^{r-1} \left(\frac{c}{l}\right)^{\gamma-1-i} H_y H_x\left(\frac{\partial^i}{\partial t^i} f(x,y,0)\right) \quad (32)$$

$$H_t H_y\ H_x\left(D_x^{\zeta} D_y^{\beta} D_t^{\gamma} f(x,y,t)\right) = \left(\frac{a}{h}\right)^{\zeta}\left(\frac{b}{k}\right)^{\gamma}\left(\frac{c}{l}\right)^{\gamma} H_t H_y\ H_x f(x,y,t) -$$

$$\begin{bmatrix} \sum_{u=0}^{m-1}\left(\frac{a}{h}\right)^{-1-u} H_t H_y\left(\frac{\partial^u}{\partial x^u} f(0,y,t)\right) - \sum_{j=0}^{n-1}\left(\frac{b}{k}\right)^{-1-j} H_t H_x\left(\frac{\partial^j}{\partial y^j} f(x,0,t)\right) \\ -\sum_{i=0}^{r-1}\left(\frac{c}{l}\right)^{-1-i} H_x H_y\left(\frac{\partial^i}{\partial t^i} f(x,y,0)\right) \\ +\sum_{i=0}^{r-1}\ \sum_{j=0}^{n-1}\ \sum_{u=0}^{m-1}\left(\frac{a}{h}\right)^{-1-u}\left(\frac{b}{k}\right)^{-1-j}\left(\frac{c}{l}\right)^{-1-i}\left(\frac{\partial^{u+j+i}}{\partial x^u y^j t^i} f(0,0,0)\right) \end{bmatrix} \quad (33)$$

*Proof.* Proof of Theorem 4.5 is similar to the proof of Theorems 4.1,4.2,4.3and 4.4 so it is omitted here.

## 5.APPLICATIONS

**(1).Solution of fractional heat equation:-**

In this part, we solve the generalized heat equation of fractional order in 2 dimensions using Shehu Transform

**Example 5.1.** Consider the equation

$$ {}_t^cD_{0^+}^{\gamma} f(x,y,t) = \frac{1}{\pi^2}\left(D_x^2 + D_y^2\right) f(x,y,t), x > 0, y > 0. t > 0. \quad (34)$$

with initial and boundary conditions

$$\begin{cases} f(0,0,t) = f(0,y,t) = f(x,0,t=0) \\ D_x(0,y,t) = \pi E_{\gamma}(-t^{\gamma}),\ \ t > 0, y > 0 \\ D_y(x,0,t) = \pi E_{\gamma}(t^{\gamma}), f(x,y,0) = \sin(\pi y)\sin(\pi x), x,y,t > 0 \end{cases} \quad (35)$$

Where $0 < \gamma \le 1$ and $D_y, D_x$ denote the corresponding partial derivative with respect to $t.y.x$ respectively.

**Solution:**

Now by using the definition of double Shehu transform on initial and boundary conditions equation(35), we get

$$H_t H_t f(0,0,t) = H_t H_y f(0,y,t) = H_t H_x\ f(x,0,t) = 0$$

$$H_t H_y D_x(0,y,t) = \frac{\pi\left(\frac{b}{k}\right)^{\gamma-1}}{\left(\frac{c}{l}\right)\left(1+\left(\frac{b}{k}\right)^{\gamma}\right)}\ ,\ \ H_t H_x D_y(x,0,t) = \frac{\pi\left(\frac{a}{h}\right)^{\gamma-1}}{\left(\frac{c}{l}\right)\left(1-\left(\frac{a}{h}\right)^{\gamma}\right)}$$

$$H_y H_x f(x,y,0) = \frac{\pi^2}{\left(\left(\frac{a}{h}\right)^2+\pi^2\right)\left(\left(\frac{b}{k}\right)^2+\pi^2\right)} \quad (36)$$

The application of triple Shehu transform on two dimension fractional heat equation and by the use of linearity property at equation (34) yields

$$H_t H_y H_x\{D_t^{\gamma} f(x,y,t)\} = \frac{1}{\pi^2}\left[H_t H_y H_x\{D_x^2 f(x,y,t)\} + H_t H_y H_x\{D_y^2 f(x,y,t)\}\right] \quad (37)$$

Using the definition of triple Shehu transform as in theorem 4.5, equation (37) reduce to the Following equation:

$$\left(\frac{c}{l}\right)^{\gamma} H_tH_yH_x\{f(x,y,t)\} = \left(\frac{c}{l}\right)^{\gamma-1} H_yH_x\{f(x,y,0)\} + \frac{1}{\pi^2}\left[\left(\frac{a}{h}\right)^2 H_tH_yH_x\{f(x,y,t)\} - H_tH_y\{f(0,y,t)\} - \left(\frac{a}{h}\right) H_tH_y\{D_x(0,y,t)\} + \left(\frac{b}{k}\right)^2 H_tH_yH_x\{f(x,y,t)\} - H_tH_x\{f(x,0,t)\} - \left(\frac{b}{k}\right) H_tH_x\{D_y(x,0,t)\}\right] \quad (38)$$

Direct substitution of equation (35) in equation (38) with some algebraic manipulation and little rearrangement yields

$$H_tH_yH_x\{f(x,y,t)\} = F[(a,b,c),(h,k,l)] = \frac{\pi^4\left(\frac{c}{l}\right)^{\gamma-1}}{\left[\left(\frac{a}{h}\right)^2+\pi^2\right]\left[\left(\frac{b}{k}\right)^2+\pi^2\right]\left[\pi^2\left(\frac{c}{l}\right)^{\gamma}-\left(\frac{a}{h}\right)^2-\left(\frac{b}{k}\right)^2\right]}$$

$$- \frac{\pi\left(\frac{a}{h}\right)\left(\frac{b}{k}\right)^{\gamma-1}}{\left(\frac{c}{l}\right)\left(1+\left(\frac{b}{k}\right)^{\gamma}\right)\left[\pi^2\left(\frac{c}{l}\right)^{\gamma}-\left(\frac{a}{h}\right)^2-\left(\frac{b}{k}\right)^2\right]} - \frac{\pi\left(\frac{b}{k}\right)\left(\frac{a}{h}\right)^{\gamma-1}}{\left(\frac{c}{l}\right)\left(1-\left(\frac{a}{h}\right)^{\gamma}\right)\left[\pi^2\left(\frac{c}{l}\right)^{\gamma}-\left(\frac{a}{h}\right)^2-\left(\frac{b}{k}\right)^2\right]} \quad (39)$$

Before to apply the inverse triple Shehu transform on equation (39), we give the following result:

$$(1+t)^{\gamma} = \sum_{u=0}^{\infty}\frac{1}{u!}\frac{\Gamma(u-\gamma)}{\Gamma(\gamma)}(-t)^u \quad (40)$$

The application of the above equation (40) in the equation (39) leads to the following equation:

$$H_tH_yH_x\{f(x,y,t)\} = F[(a,b,c),(h,k,l)]$$

$$= \sum_{u=0}^{\infty}\sum_{v=0}^{\infty}\sum_{n=0}^{\infty}\sum_{m=0}^{\infty}\frac{(-1)^{v+n+m}\,\Gamma(m-u)\pi^{-(2u+2v+2n+2)}}{m!\,\Gamma(-1)\Gamma(-1)\Gamma(-1)\Gamma(u)\left(\frac{a}{h}\right)^{2m-2u-2v}\left(\frac{b}{k}\right)^{-2m-2n}\left(\frac{c}{l}\right)^{\gamma u+1}}$$

$$-\sum_{u=0}^{\infty}\sum_{m=0}^{\infty}\sum_{r=0}^{\infty}\frac{(-1)^{v+n+m}\,\Gamma(m-u)\pi^{-(1+2u)}}{m!\,\Gamma(-1)\Gamma(-1)\Gamma(u)\left(\frac{a}{h}\right)^{2m-2u-1}\left(\frac{b}{k}\right)^{-2m-\gamma r-\gamma+1}\left(\frac{c}{l}\right)^{\gamma u+\gamma+1}}$$

$$-\sum_{u=0}^{\infty}\sum_{m-0}^{\infty}\frac{(-1)^m\,\Gamma(m-u)\pi^{-(1+2u)}}{m!\Gamma(-1)\Gamma(u)\left(\frac{a}{h}\right)^{2m-2u-\gamma+1}\left(\frac{b}{k}\right)^{-2m-1}\left(\frac{c}{l}\right)^{\gamma u+\gamma}} \quad (41)$$

Now taking the inverse triple Shehu on both sides of the equation (41), we get

$$f(x,y,t)$$

$$= \sum_{u=0}^{\infty}\sum_{v=0}^{\infty}\sum_{n=0}^{\infty}\sum_{m=0}^{\infty}\frac{(-1)^{v+n+m}\,\Gamma(m-u)\pi^{-(2u+2v+2n+2)}x^{2m-2u-2}y^{-2m-2n-1}t^{\gamma u}}{m!\,\Gamma(-1)\Gamma(-1)\Gamma(-1)\Gamma(u)\Gamma(2m-2u-2v)\Gamma(-2m-2n)\Gamma(\gamma u+1)}$$

$$-\sum_{u=0}^{\infty}\sum_{m=0}^{\infty}\sum_{r=0}^{\infty}\frac{(-1)^{m+r}\,\Gamma(m-u)\pi^{-(1+2u)}x^{2m-u-2}y^{-2m-\gamma r-\gamma}t^{\gamma u+\gamma}}{m!\,\Gamma(-1)\Gamma(-1)\Gamma(u)\Gamma(2m-2u-1)\Gamma(-2m-\gamma r-\gamma+1)\Gamma(\gamma u+\gamma+1)}$$

$$-\sum_{u=0}^{\infty}\sum_{m-0}^{\infty}\frac{(-1)^m\,\Gamma(m-u)\pi^{-(1+2u)}x^{2m-2u-\gamma}y^{-2m-2}t^{\gamma u+\gamma-1}}{m!\Gamma(-1)\Gamma(u)\Gamma(2m-2u-\gamma+1)\Gamma(-2m-1)\Gamma(\gamma u+\gamma)} \quad (42)$$

By using the definition 2.5 we can write eq.(42) as the form:

$$f(x,y,t)=\sum_{u=0}^{\infty}\quad\sum_{v=0}^{\infty}\quad\sum_{n=0}^{\infty}\frac{(-1)^{v+n}t^{\gamma u+n}x^{-2u-2}y^{-2n-1}}{\pi^{2u+2v+2n+2}}*{}'\mathrm{H}_{1,8}^{1,1}\left[-\left(\frac{x}{y}\right)^2\Big|_{(1+2u+2v,2),(1+2n,-2),(\gamma u,0)}^{(1+u,1)\atop (0,1),(2,0),(2,0),(2,0),(1-u,0)}\right]$$

$$-\sum_{u=0}^{\infty}\quad\sum_{r=0}^{\infty}\frac{t^{\gamma+\gamma u-1}x^{2u-\gamma}y^{-2}}{\pi^{1+2u}}*{}'\mathrm{H}_{1,7}^{1,1}\left[-\left(\frac{x}{y}\right)^2\Big|_{(2+2u,2),(\gamma+\gamma r,2),(-\gamma-\gamma u,0)}^{(1+u,1)\atop (0,1),(2,0),(2,0),(1-u,0)}\right]$$

$$-\sum_{u=0}^{\infty}\frac{t^{\gamma+\gamma u-1}x^{2u-\gamma}y^{-2}}{\pi^{1+2u}}*{}'\mathrm{H}_{1,6}^{1,1}\left[-\left(\frac{x}{y}\right)^2\Big|_{(2u+\gamma,2),(2,-2),(1-\gamma-\gamma u,0)}^{(1+u,1)\atop (0,1),(2,0),(1-u,0)}\right]$$

(**2**). **Solution of fractional telegraph equation** using Shehu Transform In this part, we solve The telegraph equation of fractional order in 2 dimensions using Shehu Transform

**Example 5.2** Consider the equation

$$D_t^{2\gamma}f(x,y,t)+2\alpha D_t^{\gamma}f(x,y,t)+\beta^2 f(x,y,t)=\left(D_x^2+D_y^2\right)f(x,y,t),x>0,y>0.t>0,\alpha>0,\beta>0 \quad (43)$$

with initial and boundary conditions

$$\begin{cases} f(0,y,t)=f(x,y,t)=f(x,0,t)=0 \\ D_x(0,y,t)=\pi E_\gamma(-y^\gamma),\ \ t>0,y>0 \\ D_y(x,0,t)=\pi E_\gamma(x^\gamma),f(x,y,0)=\mathrm{e}^{-y},x,y,t>0 \end{cases} \quad (44)$$

Where $0<\gamma\leq 1$ and $D_y,D_x,D_t$ denote the corresponding partial derivative with respect to $t.y.x$ respectively.

**Solution: -**

Now by using the definition of double Shehu transform on initial and boundary conditions equation (44),

we get

$$H_tH_yf(0,y,t)=H_tH_xf(x,0,t)=H_yH_x\,f(x,y,0)=0$$

$$H_tH_yD_x(0,y,t)=\frac{\pi\left(\frac{b}{k}\right)^{\gamma-1}}{\left(\frac{c}{l}\right)\left(1+\left(\frac{b}{k}\right)^{\gamma}\right)}\ ,\ \ H_tH_xD_y(x,0,t)=\frac{\pi\left(\frac{a}{h}\right)^{\gamma-1}}{\left(\frac{c}{l}\right)\left(\left(\frac{a}{h}\right)^{\gamma}-1\right)}$$

$$H_yH_xf(x,y,0)=\frac{1}{\left(\frac{a}{h}\right)\left(\frac{b}{k}+1\right)} \quad (45)$$

The application of triple Shehu transform on two dimension fractional telegraph equation and by the use of linearity property at equation (43), yields

$$H_tH_yH_x\{D_t^{2\gamma}f(x,y,t)\}+H_tH_yH_x\{2\alpha D_t^{\gamma}f(x,y,t)\}+\beta^2H_tH_yH_x\{D_t^{2\gamma}f(x,y,t)\}$$
$$=\left[H_tH_yH_x\{D_x^2f(x,y,t)\}+H_tH_yH_x\{D_y^2f(x,y,t)\}\right] \quad (46)$$

Using the definition of triple Shehu transform as in theorem 4.5, equation (46) reduce to the

Following equation:

$(1+2\alpha)\left(\frac{c}{l}\right)^{\gamma} H_t H_y H_x\{f(x,y,t)\} + \beta^2 H_t H_y H_x\{f(x,y,t)\} = (1+2\alpha)\left(\frac{c}{l}\right)^{\gamma-1} \; H_y H_x\{f(x,y,0)\} + \left(\frac{a}{h}\right)^2 H_t H_y H_x\{f(x,y,t)\} - H_t H_y\{f(0,y,t)\} - \left(\frac{a}{h}\right) H_t H_y\{D_x(0,y,t)\} + \left(\frac{b}{k}\right)^2 H_t H_y H_x\{f(x,y,t)\} - H_t H_x\{f(x,0,t)\} - \left(\frac{b}{k}\right) H_t H_x\{f(x,0,t)\}$ (47)

Direct substitution of equation (45) in equation (47) with some algebraic manipulation and little rearrangement yields

$$H_t H_y H_x\{f(x,y,t)\} = F[(a,b,c),(h,k,l)] = \frac{(1+2\alpha)\left(\frac{c}{l}\right)^{\gamma-1}}{\left(\frac{a}{h}\right)\left(\left(\frac{b}{k}\right)+1\right)\left((1+2\alpha)\left(\frac{c}{l}\right)^{\gamma}+\beta^2-\left(\frac{a}{h}\right)^2-\left(\frac{b}{k}\right)^2\right)}$$

$$- \frac{\pi\left(\frac{a}{h}\right)\left(\frac{b}{k}\right)^{\gamma-1}}{\left(\frac{c}{l}\right)\left(1+\left(\frac{b}{k}\right)^{\gamma}\right)\left((1+2\alpha)\left(\frac{c}{l}\right)^{\gamma}+\beta^2-\left(\frac{a}{h}\right)^2-\left(\frac{b}{k}\right)^2\right)} - \frac{\pi\left(\frac{b}{k}\right)\left(\frac{a}{h}\right)^{\gamma-1}}{\left(\frac{c}{l}\right)\left(\left(\frac{a}{h}\right)^{\gamma}-1\right)\left((1+2\alpha)\left(\frac{c}{l}\right)^{\gamma}+\beta^2-\left(\frac{a}{h}\right)^2-\left(\frac{b}{k}\right)^2\right)} \quad (48)$$

Application of eq.(40) in the eq.(48) leads to the following equation

$$H_t H_y H_x\{f(x,y,t)\} = F[(a,b,c),(h,k,l)]$$

$$= \sum_{p=0}^{\infty}\sum_{q=0}^{\infty}\sum_{u=0}^{\infty}\sum_{v=0}^{\infty}\sum_{n=0}^{\infty}\sum_{m=0}^{\infty}\frac{(-1)^{v+n+m}\,\Gamma(m-u)\pi^{-(2u+2v+2n+2)}\alpha^{p+q}\beta^{2q}}{m!\,\Gamma(-1)\Gamma(-1)\Gamma(-1)\Gamma(u)\left(\frac{a}{h}\right)^{2m-2u-2v}\left(\frac{b}{k}\right)^{-2m-2n}\left(\frac{c}{l}\right)^{\gamma u+1}}$$

$$-\sum_{p=0}^{\infty}\sum_{q=0}^{\infty}\sum_{u=0}^{\infty}\sum_{m=0}^{\infty}\sum_{r=0}^{\infty}\frac{(-1)^{v+n+m}\,\Gamma(m-u)\pi^{-(1+2u)}\alpha^{p+q}\beta^{2q}}{m!\,\Gamma(-1)\Gamma(-1)\Gamma(u)\left(\frac{a}{h}\right)^{2m-2u-1}\left(\frac{b}{k}\right)^{-2m-\gamma r-\gamma+1}\left(\frac{c}{l}\right)^{\gamma u+\gamma+1}}$$

$$-\sum_{p=0}^{\infty}\sum_{q=0}^{\infty}\sum_{u=0}^{\infty}\sum_{m-0}^{\infty}\frac{(-1)^{m}\,\Gamma(m-u)\pi^{-(1+2u)}\alpha^{p+q}\beta^{2q}}{m!\Gamma(-1)\Gamma(u)\left(\frac{a}{h}\right)^{2m-2u-\gamma+1}\left(\frac{b}{k}\right)^{-2m-1}\left(\frac{c}{l}\right)^{\gamma u+\gamma}} \quad (49)$$

Now, taking the inverse triple Shehu on both sides of the equation (49) we get

$$f(x,y,t) = \sum_{p=0}^{\infty}\sum_{q=0}^{\infty}\sum_{u=0}^{\infty}\sum_{v=0}^{\infty}\sum_{n=0}^{\infty}\sum_{m=0}^{\infty}\frac{(-1)^{v+n+m}\,\Gamma(m-u)\Gamma(p+q)}{m!\,\Gamma(-1)\Gamma(-1)\Gamma(-1)\Gamma(u)}$$

$$*\frac{x^{2m-2u-2v-1}y^{-2m-2n-1}t^{\gamma u}\alpha^{p+q}\beta^{2q}}{\Gamma(2m-2u-2v)\Gamma(-2m-2n)\Gamma(\gamma u+1)\pi^{(2u+2v+2n+2)}}$$

$$-\sum_{p=0}^{\infty}\sum_{q=0}^{\infty}\sum_{u=0}^{\infty}\sum_{m=0}^{\infty}\sum_{r=0}^{\infty}\frac{(-1)^{m+r}\,\Gamma(m-u)\Gamma(p+q)}{m!\,\Gamma(-1)\Gamma(-1)\Gamma(u)}$$

$$*\frac{x^{2m-u-2}y^{-2m-\gamma r-\gamma}t^{\gamma u+\gamma}\alpha^{p+q}\beta^{2q}}{\Gamma(2m-2u-1)\Gamma(-2m-\gamma r-\gamma+1)\Gamma(\gamma u+\gamma+1)\pi^{(1+2u)}}$$

$$-\sum_{p=0}^{\infty}\sum_{q=0}^{\infty}\sum_{u=0}^{\infty}\sum_{m-0}^{\infty}\frac{(-1)^{m}\,\Gamma(m-u)\Gamma(p+q)x^{2m-2u-\gamma}y^{-2m-2}t^{\gamma u+\gamma-1}\alpha^{p+q}\beta^{2q}}{m!\Gamma(-1)\Gamma(u)\Gamma(2m-2u-\gamma+1)\Gamma(-2m-1)\Gamma(\gamma u+\gamma)\pi^{(1+2u)}} \quad (50)$$

By using definition 2.5 we can write eq.(50) as the form:

$$f(x,y,t)=\sum_{p=0}^{\infty}\sum_{q=0}^{\infty}\sum_{u=0}^{\infty}\sum_{v=0}^{\infty}\sum_{n=0}^{\infty}\frac{(-1)^{v+n}t^{\gamma u}x^{-2u-2}y^{-2n-1}\alpha^{p+q}\beta^{2q}}{\pi^{2u+2v+2n+2}}$$

$$*\,'\mathrm{H}_{1,8}^{2,1}\left[-\left(\frac{x}{y}\right)^{2}\middle|\begin{matrix}(1+u,1),(1-p-q,0)\\(0,1),(2,0),(2,0),(2,0),(1-u,0)\\(1+2u+2v,2),(1+2n,-2),(\gamma u,0)\end{matrix}\right]$$

$$f(x,y,t)=\sum_{p=0}^{\infty}\sum_{q=0}^{\infty}\sum_{u=0}^{\infty}\sum_{r=0}^{\infty}\frac{(-1)^{r}t^{\gamma+\gamma u}x^{-u-2}y^{-\gamma-\gamma r}\alpha^{p+q}\beta^{2q}}{\pi^{1+2u}}$$

$$*\,'\mathrm{H}_{1,7}^{2,1}\left[-\left(\frac{x}{y}\right)^{2}\middle|\begin{matrix}(1+u,1),1-p-q,0)\\(0,1),(2,0),(2,0),(1-u,0)\\(2+2u,2),(\gamma+\gamma r,2),(-\gamma-\gamma u,0)\end{matrix}\right]$$

$$f(x,y,t)=\sum_{p=0}^{\infty}\sum_{q=0}^{\infty}\sum_{u=0}^{\infty}\frac{t^{\gamma+\gamma u-1}x^{2u-\gamma}y^{-2}\alpha^{p+q}\beta^{2q}}{\pi^{1+2u}}$$

$$*\,'\mathrm{H}_{1,6}^{2,1}\left[-\left(\frac{x}{y}\right)^{2}\middle|\begin{matrix}(1+u,1),1-p-q,0)\\(0,1),(2,0),(1-u,0)\\(2u+\gamma,2),(2,-2),(1-\gamma-\gamma u,0)\end{matrix}\right] \qquad (51)$$

## 6.CONCLUSIONS

The triple Shehu transform formulas for fractional partial integrals and derivatives in the sense of Caputo fractional derivative have been establishedin this paper. After that, a few theorems about this new transform have been proved. In addition, the double Shehu transform was used to solve several linear partial fractional differential equations. It's worth mentioning that(TRHT) can be used in combination with other approaches to solve non- linear (PDEs) in applied mathematics, engineering, and applied physics, allof which will be addressed in future articles.

Wagdi F. S. Ahmed
School of Mathematical Sciences, Swami Ramanand Teerth
Marathwada University, Nanded-431606, India.
*e-mail: wagdialakel@gmail.com*
D. D. Pawar
School of Mathematical Sciences, Swami Ramanand Teerth
Marathwada University, Nanded-431606, India.
e-mail: dypawar@yahoo.com